\documentclass[11pt,reqno]{amsart}

\usepackage{amssymb}
\usepackage{latexsym}
\usepackage{amsmath}

\def\wt{\widetilde}

\def\ov{\overline}

\def\up{\upharpoonright}
\def\cH{\mathcal H}
\def\cD{\mathcal D}
\def\cR{\mathcal R}

\def \gH{\mathfrak H}

\def \gN{\mathfrak N}
\def \bC{\mathbb C}

\def \bR{\mathbb R}
\def \bN{\mathbb N}

\def\cB{\mathcal B}

\def \cd {\cdot}

\def \a{\alpha}
\def \b{\beta}
\def \l{\lambda}
\def \s{\sigma}
\def \t{\theta}

\def \f{\varphi}
\def \G{\Gamma}

\def \Ker{\text {Ker}}

\def \im{\text{Im}}
  \def\ri {\wt\rho_I (A)}

\def\bt{\{\cH,\Gamma_0,\Gamma_1\}}

\newtheorem{theorem}{Theorem}[section]
\newtheorem{proposition}[theorem]{Proposition}
\newtheorem{corollary}[theorem]{Corollary}
\newtheorem{lemma}[theorem]{Lemma}

\theoremstyle{definition}
\newtheorem{definition}[theorem]{Definition}

\theoremstyle{remark}
\newtheorem{remark}[theorem]{Remark}
\numberwithin{equation}{section}

\begin{document}


\title [Symmetric operators]
{Symmetric operators with real defect subspaces of the maximal dimension.
Applications to differential operators}

\author{Vadim Mogilevskii}
\address{Department of Math. Analysis,
Lugans'k National   University, 2 Oboronna str., Lugans'k  91011, Ukraine}

\email{vim@mail.dsip.net}



\begin{abstract}
Let $\gH$ be a Hilbert space and let $A$ be a simple symmetric operator in
$\gH$ with equal  deficiency indices $d:=n_\pm(A)<\infty$. We show that if, for
all $\l$ in an open interval $I\subset\bR$, the dimension of defect subspaces
$\gN_\l(A)(=\Ker (A^*-\l))$ coincides with $d$, then every self-adjoint
extension $\wt A\supset A$ has no continuous spectrum in $I$ and the point
spectrum of $\wt A$ is nowhere dense in $I$. Application of this statement to
differential operators makes it possible to generalize the known results by
Weidmann to the case of an ordinary differential expression with both singular
endpoints and arbitrary equal deficiency indices of the minimal operator.
Moreover, we show in the paper, that an old conjecture by Hartman and Wintner
on the spectrum of a self-adjoint Sturm - Liouville operator is not valid.
\end{abstract}


\maketitle

\section{Introduction}
 Let $\gH$ be a Hilbert space, let $A$ be a simple symmetric densely
defined operator in $\gH$ with equal and finite deficiency indices
$d=n_\pm(A)<\infty$ and let $\gN_z(A)=\Ker (A^*-z), \; z\in\bC $ be a defect
subspace of $A$. As is known \cite{Nai} $\dim \gN_\l(A)\leq d$ for all
$\l\in\bR$ and $\dim \gN_\l(A)= d$ if the range of $A-\l$ is closed, i.e., if
$\l$ belongs to the set $\hat\rho (A)$ of all regular type points of $A$ (note
that $\Ker (A-\l)=\{0\}$, since the operator $A$ is simple). Moreover, if
$I=(\mu_1,\mu_2)$ is an interval such that $I\subset \hat\rho (A)$, then for
any self-adjoint extension $\wt A\supset A$ the spectrum $\s (\wt A)$ in $I$
consists of isolated eigenvalues  of $\wt A $ with finite multiplicity (the
discrete spectrum $\s_d(\wt A)$). In this connection it seems to be rather
interesting to find out if the situation is the same for the weaker condition
\begin {equation}\label{0.1}
\dim \gN_\l(A)= d, \quad \l\in I.
\end{equation}
It turns out that the answer is negative. More precisely, we show in the paper
(see Proposition \ref{pr8}) that for any interval $I$ there is an operator $A$
such that \eqref{0.1} is satisfied and for any (equivalently for some)
self-adjoint extension $\wt A\supset A$ the set of all points $\l\in I$
belonging to the essential spectrum $\s_e(\wt A)(=\s (\wt A)\setminus \s_d (\wt
A))$ is infinite. At the same time the spectrum of such an extension $\wt A$ is
"small" enough. Namely, in the main theorem of the paper we prove that under
the condition \eqref{0.1} the following statement (s) is valid for any
self-adjoint extension $\wt A\supset A$:

 \hskip 3 mm (s)  the set $\s (\wt A)\cap I$ is
nowhere dense in $I$ and coincides with the closure of the set $\s_p (\wt
A)\cap I$, where $\s_p (\wt A)$ is the set of all eigenvalues of $\wt A$ (the
point spectrum).

Our considerations  are substantially inspired by  the book \cite{Wei} and the
recent paper \cite{Zet08} in Journal of  Funct. Anal., where similar results
were obtained for differential operators. Namely, let $L_0$ be the minimal
symmetric operator generated by a formally self-adjoint differential expression
$l[y]$ of an even order $2n$ on an interval $(a,b), -\infty\leq a<b\leq\infty$
(see \eqref{29}). For the operator $L_0$ satisfying \eqref{0.1} the validity of
the statement (s) for any extension $\wt A=\wt A^* \supset L_0$ was proved by
Weidmann \cite{Wei} under the assumptions, that $a$ is a regular endpoint for
the expression $l[y]$ and $L_0$ has minimal deficiency indices
$d(=n_\pm(L_0))=n$. Moreover, it was shown in \cite{Zet08} that in the case  of
the regular endpoint $a$ and  an arbitrary defects $d$ the statement (s) holds
for some self-adjoint extension $\wt A\supset L_0$ defined by separated
boundary conditions.

In the present paper we generalize the Weidmann's result to the case of
arbitrary (regular or singular) endpoints $a$ and $b$ and arbitrary equal
deficiency indices $d=n_\pm (L_0)$. More precisely, let $L_{a0}$ and $L_{b0}$
be minimal operators for the expression $l[y]$ on intervals $(a,c)$ and $(c,b)$
respectively (with some $c\in (a,b)$), let $n_+(L_{a0})=n_-(L_{a0})=:d_a, \;
n_+(L_{b0})=n_-(L_{b0})=:d_b$ and let for some interval $I=(\mu_1,
\mu_1)\subset \bR$
\begin {equation*}
\dim\,\gN_\l (L_{a0})=d_a, \quad \dim\,\gN_\l (L_{b0})=d_b, \quad \l\in I.
\end{equation*}
We show in Theorem \ref{th9} that under such assumptions the statement (s)
holds for any self-adjoint extension $\wt A\supset L_0$.

In the paper \cite{HarWin49} Hartman and Wintner  suggested that for the second
order, i.e. Sturm - Liouville,  operator $L_0$ on the semiaxis $[0,\infty)$
with
\begin {equation}\label{0.2}
\dim \gN_\l(L_0)=1(= d), \quad \l\in I=(\mu_1,\mu_2).
\end{equation}
the statement (s) can be strengthened to  "the spectrum of any self-adjoint
extension $\wt A\supset L_0$ is discrete in $I$" (similar conjecture for the
operator $L_0$ of an arbitrary order $2n$ is contained in \cite{Wei,Zet08}). We
prove in the paper, that this conjecture is not valid. More precisely we show
that for any finite interval $I=(\mu_1,\mu_2)$ there exists a Sturm - Liouville
operator $L_0$ such that \eqref{0.2} holds and for any self-adjoint extension
$\wt A\supset L_0$ the set $\s_e (\wt A)\cap I$ is infinite (see Proposition
\ref{pr11}).

In conclusion note that our approach is based on the concepts of a boundary
triplet for $A^*$ and the corresponding abstract Weyl function, which has
become a convenient tool in the extension theory of symmetric operators and its
applications (see \cite{GorGor,DM91,Mal92,DM95,Mog09.1} and references
therein). Such an approach enabled us to obtain the above results without
 complicated construction of the self-adjoint extension $\wt
A\supset L_0$ with the desired properties of the spectrum $\s (\wt A)$ (cf.
\cite{Zet08}).
\section{Preliminaries}
In the sequel we use the following notations: $\gH$, $\cH$ denote separable
Hilbert spaces; $[\cH_1,\cH_2]$  is the set of all bounded linear operators
defined on $\cH_1$ with values in $\cH_2$; $[\cH]:=[\cH,\cH]$; $\bC_+\,(\bC_-)$
is the upper (lower) half-plain of the complex plain. Moreover, for a (not
necessarily bounded) operator $T$ from $\cH_1$ to $\cH_2$ we denote by
$\cD(T),\,\cR (T)$ and $\text {Ker}T$ the domain,  range and the kernel of $T$
respectively.

For a closed operator $T$ in $\gH$ we denote by $\hat\rho (T)=\{\l \in \bC:\Ker
(T-\l)=\{0\}, \; \ov{\cR (T-\l)}=\cR (T-\l)\}$ and $\rho (T)=\{\l \in \hat\rho
(T) :\cR (T-\l)=\gH\}$ the set of regular type points and the resolvent set of
$T$ respectively.

Let $\cH$ be a finite dimensional Hilbert space. Recall that a holomorphic
operator function $\Phi (\cd):\bC_+\cup\bC_-\to [\cH]$ is called a Nevanlinna
function (and is referred to the class $R[\cH]$) if $\im\, z\cd \im \Phi
(z)\geq 0 $ and $\Phi ^*(z)= \Phi (\ov z), \; z\in\bC_+\cup\bC_-$. According to
\cite{KacKre,Br} a function $\Phi (\cd):\bC_+\cup\bC_-\to [\cH]$ belongs to the
class $R[\cH]$ if and only if it admits the integral representation
\begin {equation}\label{1}
\Phi (z)= C_0 +z\,C_1 +\int_\bR \left(\frac 1 {t-z} - \frac t
{1+t^2}\right)dF_\Phi (t),
\end{equation}
where $C_0,C_1\in [\cH], \; C_0=C_0^*, \; C_1\geq 0$ and $F_\Phi (\cd):\cB_b\to
[\cH]$ is an operator valued measure defined on the ring $\cB_b$ of all bounded
Borel sets in $\bR$ and such that
\begin {equation}\label{2}
\int_\bR \frac 1 {t^2+1}dF_\Phi (t) \in [\cH].
\end{equation}
The operator valued measure $F_\Phi (\cd)$ in \eqref{1} is called the spectral
measure of the function $\Phi (\cd)\in R [\cH]$.

The following lemma is well known.
\begin{lemma}\label{lem1}
Let $\Phi (\cd)\in R[\cH]$ and let $F(\cd)=F_\Phi(\cd)$ be the corresponding
spectral measure. Then for each $\l\in\bR$ the following relations are
equivalent:
\begin{gather}
\text{(i)} \quad \lim_{y\to 0}  1/ y \;\im\, (\Phi (\l+iy)h,h)<\infty, \quad
h\in\cH; \nonumber\\
\text{(ii)}\quad \int_\bR \frac {d(F(t)h,h)} {(t-\l)^2}< \infty,\quad
h\in\cH;\qquad\qquad \label{3}
\end{gather}
If the relation (i) (or, equivalently, (ii)) is satisfied, then there exists
the limit
\begin {equation*}
M(\l+i0):=\lim_{y\to 0} M(\l+iy)
\end{equation*}
and $\im\, M(\l+i0)=0$.
 \end{lemma}
Let $A$ be a closed densely defined symmetric operator in $\gH$ and let $A^*$
be the adjoint operator. For each $z\in\bC$ denote by
\begin {equation*}
\gN_z(A):=\Ker (A^*-z)(=\gH\ominus \cR (A-\ov z))
\end{equation*}
the defect subspace of $A$ and let $n_\pm (A)=\dim \gN_z(A)\; (z\in\bC_\pm)$ be
the deficiency indices of $A$.
\begin{definition}\  \cite{GorGor}
A triplet $\Pi=\bt$ consisting of an auxiliary Hilbert space $\cH$ and linear
mappings $\G_j:\cD (A^*)\to\cH, \; j\in \{0,1\}$ is called a boundary triplet
for $A^*$ if the mapping $\G=(\G_0\;\;\G_1)^\top: \cD (A^*)\to\cH\oplus\cH$ is
surjective and the following abstract Green's identity holds
\begin {equation*}
(A^*f,g)-(f,A^*g)=(\G_1 f,\G_0 g) - (\G_0 f,\G_1 g), \quad f,g\in \cD (A^*).
\end{equation*}
\end{definition}
The following Proposition was proved in \cite{DM92}.
\begin{proposition}\label{pr3}
Let $\Pi=\bt$ be a boundary triplet for $A^*$. Then $n_+(A)=n_-(A)=\dim\cH$ and
the equalities
\begin {equation}\label{4}
\cD (A_0)=\Ker \G_0=\{f\in\cD (A^*):\G_0 f=0\}, \qquad A_0=A^* \up \cD (A_0)
\end{equation}
define a self-adjoint extension $A_0\supset A$.

Conversely, let $A$ be a symmetric operator in $\gH$ with $n_+(A)=n_-(A)$ and
let  $\wt A $ be a self-adjoint extension of $A$. Then there exists a boundary
triplet $\Pi=\bt$ for $A^* $ such that $\wt A=A_0(=A^*\up \Ker \G_0)$.
\end{proposition}
It turns out that for any $z\in \rho (A_0)$ the operator $\G_0\up \gN_z(A)$
isomorphically maps $\gN_z(A)$ onto $\cH$. This enables one to introduce the
following definition.
\begin{definition}\ \cite{DM91}
The operator function $M(\cd):\rho (A_0)\to [\cH]$ defined by
\begin {equation*}
\G_1\up \gN_z(A)= M(z)\G_0\up \gN_z(A), \quad z\in\rho (A_0)
\end{equation*}
is called the  Weyl function corresponding to the boundary triplet $\bt$.
\end{definition}
As was shown in \cite{DM91} the Weyl function $M(\cd)$ belongs to the class
$R[\cH]$ and $0\in\rho (\im M(z)), \; z\in\bC_+\cup\bC_-$.

\section{Symmetric operators with real defect subspaces of the maximal dimension}
In the sequel we denote by $A$ a simple symmetric densely defined operator in
$\gH$ with equal deficiency indices $d=n_\pm (A)<\infty$. Since the operator
$A$ is simple, it follows that $\Ker (A-\l)=\{0\}$ and, consequently, $\dim
\gN_\l (A)\leq d$ for all $\l\in\bR$. We denote by $\wt\rho (A)$ the set of all
$\l\in\bR$ such that $\dim \gN_\l (A)=d$.
\begin{proposition}\label{pr5}
Assume that $A$ is a simple symmetric operator in $\gH$ with $d=n_\pm
(A)<\infty$, $\Pi=\bt$ is a boundary triplet for $A^*$, $A_0$ is the
self-adjoint extension \eqref{4} and $M(\cd)$ is the corresponding Weyl
function. Then a real point $\l$ belongs to $\wt\rho (A)$ and $\Ker
(A_0-\l)=\{0\}$ if and only if
\begin {equation}\label{9}
\lim_{y\to 0} 1/y \,\im (M(\l+iy)h,h)<\infty, \quad h\in \cH.
\end{equation}
\end{proposition}
\begin{proof}
For a point $\l\in\bR$ denote by $\cH_\l$ the subspace in $\cH$ given by
$\cH_\l=\G_0\gN_\l (A)$. It follows from \eqref{4} that
\begin {equation*}
\Ker (\G_0\up \gN_\l(A))=\cD (A_0)\cap\gN_\l(A)=\Ker (A_0-\l)
\end{equation*}
and, consequently,
\begin {equation}\label{10}
\dim \gN_\l(A)=\dim \Ker (A_0-\l)+\dim \cH_\l.
\end{equation}
Since $\dim\gN_\l(A)\leq d$, the equality \eqref{10} yields the equivalences
\begin {equation}\label{11}
(\Ker (A_0-\l)=\{0\} \;\;\text {and} \;\; \dim \gN_\l(A)=d)\iff
\dim\cH_\l=d\iff \cH_\l=\cH.
\end{equation}
Moreover according to \cite{Mal92} for any $h\in\cH$ the following equivalence
holds
\begin {equation}\label{12}
h\in \cH_\l\iff \lim_{y\to 0} 1/y \,\im (M(\l+iy)h,h)<\infty.
\end{equation}
In view of \eqref{12} the equality $\cH_\l=\cH$ is equivalent to the condition
\eqref{9}, which together with \eqref{11} gives the desired statement.
\end{proof}
\begin{remark}
It is easily seen that for all $\l\in\bR$ there exists a self-adjoint extension
$A_0\supset A$ with $\Ker (A_0-\l)=\{0\}$. Moreover by Proposition \ref{pr3}
there exists a boundary triplet $\Pi=\bt$ for $A^*$ such that $A_0$ is defined
by \eqref{4}. This implies that Proposition \ref{pr5} provides actually the
criterium (in terms of the Weyl function) for a point $\l\in\bR$ belongs to
$\wt\rho (A)$.
\end{remark}
\begin{lemma}\label{lem6}
Assume that $\dim\cH<\infty$ and $F(\cd):\cB_b\to [\cH]$ is an operator valued
measure satisfying \eqref{2} and the relation
\begin {equation*}
\lim_{[\alpha,\beta)\to \bR} (F([\alpha,\beta))h,h)=\infty, \quad h\in\cH.
\end{equation*}
Moreover, let $L_2(F,\cH)$ be the Hilbert space of vector functions
$f(\cd):\bR\to \cH$ such that
\begin {equation*}
||f||^2_{L_2(F,\cH)}:=\int_\bR (d F(t)f(t),f(t))<\infty
\end{equation*}
(see \cite{Kac50,DunSch}) and let $\wt A_F$ be the self-adjoint multiplication
operator  in $L_2(F,\cH)$ given by
\begin {equation}\label{12.1}
\cD (\wt A_F)= \{f\in L_2(F,\cH): \, t f(t)\in L_2(F,\cH)\}, \quad (\wt A_F
f)(t)=tf(t).
\end{equation}
Then: 1) the equalities
\begin {equation}\label{12.2}
\cD (A_F)=\left \{f\in \cD (\wt A_F): \int_\bR d F(t)f(t)=0\right\}, \quad (A_F
f)(t)=tf(t)
\end{equation}
define a simple symmetric densely defined operator in $L_2(F,\cH)$ such that
$n_\pm (A_F)=\dim \cH$ and $A_F\subset \wt A_F$;

2) for each point $\l\in\bR$ with $F(\{\l\})=0(\Leftrightarrow \Ker (\wt
A_F-\l)=\{0\})$ the inclusion $\l\in \wt\rho (A_F)$ is equivalent to the
relation \eqref{3}.
\end{lemma}
\begin{proof}
The statement 1) was proved in \cite{DM95}.

2) Let the function $M(\cd)\in R[\cH]$ be given by \eqref{1} with $C_1=0$ and
$F_M(\cd)=F(\cd)$. Then according to \cite{DM95} there exists a boundary
triplet $\Pi_0=\bt$ for $A_F^*$ such that $A_0(=A^*\up \Ker \G_0)=\wt  A_F$ and
the corresponding Weyl function coincides with $M(\cd)$. Applying now
Proposition \ref{pr5} to the triplet $\Pi_0$ and taking Lemma \ref{lem1} into
account one obtains the desired statement.
\end{proof}
For a given operator $A$ and an interval $I=(\mu_1,\mu_2), \; -\infty\leq\mu_1
<\mu_2\leq \infty,$ we denote by $\ri=\wt\rho (A)\cap I$ the set of all points
$\l\in I$ with $\dim\gN_\l (A)=d(=n_\pm (A))$ and let $\hat\rho_I (A)=\hat \rho
(A)\cap I $ be the set of all regular type points of $A$ belonging to $I$.
Since $\dim\gN_\l (A)=d$ for all $\l\in \hat\rho_I (A)$, the inclusion
$\hat\rho_I (A)\subset \ri$ is valid. Moreover, the set $\ri\setminus
\hat\rho_I (A)$ consists of all points $\l\in I$ such that $ \dim \gN_\l (A)=d$
and the range $\cR (A-\l)$ is not closed.

As is known \cite{RS} the spectrum $\s  (T)$ of a self-adjoint operator $T$
admits the representation
\begin {equation}\label{13}
\s (T)=\ov{\s_p (T)} \cup \s_c(T), \quad   \s_c(T)=\s_{ac}(T)\cup \s_{sc}(T),
\end{equation}
where $\s_p (T)=\{\l\in\bR: \Ker (T-\l)\neq \{0\}\}$ is the point spectrum and
$\s_c(T), \; \s_{ac}(T)$ and $\s_{sc}(T)$ are continuous, absolutely continuous
and singular continuous parts of $\s(T)$ respectively. Recall that the
continuous spectrum $\s_c(T)$ is defined as the spectrum of the self-adjoint
operator $T_c=T\up \gH_c$,  where $\gH_c:=\gH\ominus \text{span} \{Ker
(T-\l):\l\in\s_p(T)\}$ is  the subspace reducing the operator $T$.

 Another basic partition of the spectrum is in terms of the discrete spectrum
$\s_d(T)$ and the essential spectrum $\s_e(T)$. Namely, $\s_d(T)$ is the set of
all isolated eigenvalues of $T$ with finite multiplicity and
$\s_e(T)=\s(T)\setminus \s_d(T)$. It is clear that $\s_c(T)\subset \s_e(T)$.
Moreover, the following lemma is well known.
\begin{lemma}\label{lem6a}
Let $A$ be a simple symmetric operator with $d=n_\pm (A)<\infty$. Then all
self-adjoint extensions $\wt A\supset A$ have the same essential spectrum
\begin {equation}\label{13a}
\s_e(\wt A)=\bR\setminus \hat\rho (A).
\end{equation}
\end{lemma}
Recall also that a set $X\subset (\mu_1,\mu_2)$  is called nowhere dense in
$(\mu_1,\mu_2)$ if for any interval $(\mu_1',\mu_2')\subset (\mu_1,\mu_2) $
there exists an interval $(\mu_1'',\mu_2'')\subset (\mu_1',\mu_2')$ such that
$X\cap (\mu_1'',\mu_2'')=\emptyset$.

Now we are ready to prove the main theorem of the paper.

\begin{theorem}\label{th7}
Assume that $A$ is a simple symmetric densely defined operator in $\gH$ with
equal deficiency indices $d=n_\pm (A)<\infty$ and $I=(\mu_1,\mu_2), \;
-\infty\leq\mu_1 <\mu_2\leq \infty,$ is an interval such that the set
$I\setminus \ri$ is at most countable. Then:

1) for each self-adjoint extension $\wt A\supset A$ the intersection $\s_c (\wt
A)\cap I$  is empty and the set $\s (\wt A)\cap I$ is nowhere dense in $I$;

2) the set $I\setminus \hat \rho_I(A)$ is nowhere dense in $I$.
\end{theorem}
\begin{proof}
1) Let $\wt A$ be a self-adjoint extension of $A$ and let $\Pi=\bt$ be a
boundary triplet for $A^*$ with $\wt A=A_0(=A^*\up \Ker \G_0)$ (such a triplet
exists in view of Proposition \ref{pr3}). Moreover, let $M(\cd)\in R[\cH] $ be
the corresponding Weyl function.

Next assume that
\begin {equation*}
X_p=\{\l_k\}(=\s_p (\wt A)\cap I)
\end{equation*}
is the (at most countable) set of all eigenvalues of $\wt A$ belonging to $I$
and let $X_1:=\ri\setminus X_p, \; X_2:= (I\setminus \ri)\setminus X_p $, so
that
\begin {equation}\label{14}
I=X_p\cup X_1\cup X_2, \quad X_p\cap X_1=X_p\cap X_2=X_1\cap X_2=\emptyset.
\end{equation}
Then $X_p\cup X_2$ is an at most countable subset   in $I$ and by Proposition
\ref{pr5} the Weyl function $M(\cd)$ satisfies the relation \eqref{9} for all
$\l\in X_1$. This and Lemma \ref{lem1} yield the following statement
($\text{s}_1$) :

\hskip 5mm ($\text{s}_1$) there exists a subset $X_1\subset I$ such that: (i)
$I\setminus X_1$ is an at most countable set; (ii) for all $\l\in X_1$ the
limit $M(\l+i0):=\lim\limits_{y\to 0} M(\l+iy)$ exists and $\im\, M(\l+i0)=0$.

According to \cite[Theorem 4.3]{MalNei02} the statement ($\text{s}_1$) implies
that $\s_{sc}(\wt A)\cap I=\emptyset, \; \s_{ac}(\wt A)\cap I=\emptyset$ and,
consequently,
\begin {equation}\label{15}
\s_c(\wt A)\cap I=\emptyset.
\end{equation}
Next assume that $E(\cd)$ and $F(\cd)=F_M(\cd)$ are spectral measures of the
operator $\wt A(=A_0)$ and the Weyl function $M(\cd)$ respectively. According
to \cite[Lemma 3.2]{MalNei02} the measures $E(\cd)$ and $F(\cd)$ are
equivalent. Moreover, by \eqref{15}
\begin {equation}\label{16}
\s(\wt A)\cap I=\ov X_p,
\end{equation}
which implies that the measure $E(\cd)$ is discrete on $I$ and hence so is the
measure $F(\cd)$. Combining this statement with Proposition \ref{pr5} and Lemma
\ref{lem1} one obtains
\begin {equation}\label{17}
\int_\bR \frac {d(F(t)h,h)} {(t-\l)^2}=\sum_k  {(F_k h,h)}/ {(\l_k-\l)^2} <
\infty,\quad \l\in X_1,\;\; h\in\cH,
\end{equation}
where $F_k=F(\{\l_k\})\in [\cH]$. Let $\{e_j\}_1^d$ be an orthonormal basis in
$\cH$ and let $c_k=\sum\limits_j (F_k e_j, e_j)$. Since $F_k\neq 0$ and
$F_k\geq 0$, it follows that $c_k>0$ and the relation \eqref{17} yields
\begin {equation}\label{18}
\sum_k  c_k / {(\l_k-\l)^2} < \infty,\quad \l\in X_1.
\end{equation}
Thus the following statement ($\text{s}_2$) is proved:

\hskip 5mm ($\text{s}_2$) there exists a decomposition \eqref{14} of the
interval $I$ and a sequence of positive numbers $\{c_k\}$ such that
$X_p=\{\l_k\}$ and $X_2$ are at most countable sets and for all $\l\in X_1$ the
relation \eqref{18} holds.

By using the statement ($\text{s}_2$) one can prove in the same way as it was
done \cite[Theorem11.7]{Wei} that the set $X_p$ is nowhere dense in $I$. This
and the equality \eqref{16} imply that the set $\s (\wt A)\cap I$ is nowhere
dense in $I$ as well.

The statement 2) follows from the obvious inclusion $(I\setminus \hat
\rho_I(A))\subset \s (\wt A)\cap I$ and the statement 1).
\end{proof}
It turns out that the relation $\s_c (\wt A)\cap I=\emptyset$ in the statement
1) of Theorem \ref{th7} can not be replaced with the stronger one $\s_e (\wt
A)\cap I=\emptyset$. More precisely the following proposition holds.
\begin{proposition}\label{pr8}
For any interval $I=(\mu_1,\mu_2), \; -\infty\leq\mu_1 <\mu_2\leq \infty,$ and
for any $d\in\bN$ there exist a Hilbert space $\gH$ and a simple symmetric
operator $A$ in $\gH$ such that $n_\pm(A)=d, \; \ri=I$ and for any self-adjoint
extension $\wt A\supset A$ the interval $I$ contains infinitely many points of
$\s_e (\wt A)$. In view of \eqref{13a} the last statement implies that the set
$\ri\setminus \hat \rho_I (A) $ is infinite.
\end{proposition}
\begin{proof}
First assume that $d=1$ and consider the following two alternative cases:

\hskip 5mm (i) $\mu_2<\infty$, i.e., the interval $I$ is bounded from above.

Let $\{\l_k\}_0^\infty$ be a strictly increasing sequence of the points
$\l_k\in I $ such that $\lim\limits_{k\to\infty}\l_k=\mu_2$ and let
$\{\l_{jk}\}_{j,k=1}^\infty$ be a sequence of the points $\l_{jk}\in (\l_{k-1},
\l_k)$ such that $\l_{jk}<\l_{j+1,k}, \; \l_k-\l_{jk}<1$ and
$\lim\limits_{j\to\infty}\l_{jk}=\l_k, \; j,k \in\bN$. Consider also two
sequences of positive numbers $\{s_k\}_1^\infty$ and $\{u_{jk}\}_
{j,k=1}^\infty$ such that
\begin {equation*}
\sum_{k=1}^\infty s_k<\infty \;\;\; \text{and} \;\;\; \sum_{j=1}^\infty
u_{jk}=s_k
\end{equation*}
and let $F_{jk}=u_{jk} (\l_k-\l_{jk})^2$. Since $F_{jk}\leq u_{jk}$, it follows
that
\begin {equation*}
\sum_k \sum_j F_{j k}\leq \sum_k \sum_j u_{j k}=\sum_k  s_k <\infty.
\end{equation*}
This enables one to introduce the scalar discrete measure $F'(\cd)$ on Borel
sets $B\subset I$ by
\begin {equation}\label{20}
F'(\{\l_{jk}\})=F_{jk}, \quad F'(B)=\sum_{\l_{jk}\in B} F_{jk}.
\end{equation}
Assume also that $F''(\cd)$ is a scalar measure on bounded Borel sets $B\subset
\bR $ such that $F''(\bR\setminus I) =\infty$ and $\int\limits_{\bR\setminus I}
(t^2+1)^{-1} dF''(t)<\infty$ (for example one can take as $F''(\cd)$ the
standard Lebesgue measure on the line). Then the equality
\begin {equation}\label{21}
F(B)=F'(B\cap I)+F''(B\cap (\bR\setminus I))
\end{equation}
defines the scalar measure $F(\cd)$ on bounded Borel sets $B\subset \bR$ such
that $F(\bR)=\infty$ and
\begin {equation}\label{22}
\int_\bR\frac {d F(t)} {t^2+1} <\infty.
\end{equation}
Next we show that
\begin {equation}\label{23}
\int_\bR \frac {d F(t)} {(t-\l_k)^2} <\infty
\end{equation}
for each point $\l_k$. Let $(\a_k,\b_k)$ be an interval such that
$\l_{k-1}<\a_k<\l_{1,k}$ and $\l_{k}<\b_k<\l_{1,k+1}$. Then
\begin {equation}\label{24}
\int\limits_{[\a_k,\b_k)}\frac {d F(t)} {(t-\l_k)^2}=\sum_ {j=1}^\infty\frac {
F_{jk}} {(\l_{jk}-\l_k)^2}=\sum _{j=1}^\infty u_{jk}=s_k<\infty
\end{equation}
and in view of \eqref{22} one has
\begin {equation}\label{25}
\int\limits_{\bR\setminus [\a_k,\b_k)}\frac {d F(t)} {(t-\l_k)^2}<\infty.
\end{equation}
Combining \eqref{24} and \eqref{25} we arrive at \eqref{23}.

\hskip 5mm (ii) $\mu_2=\infty$, i.e., the interval $I$ is unbounded from above.

Assume also without loss of generality that $\mu_1\leq 1$. In this case we put
$\l_k=k, \; \l_{jk}=k+ \tfrac 1 {j+1}, \; k,j \in \bN$ and let
$\{F_j\}_1^\infty$ be a sequence of positive numbers such that
$\sum\limits_{j=1}^\infty (j+1)^2 F_j<\infty$. Then
$s:=\sum\limits_{j=1}^\infty  F_j<\infty$ and the equalities
\begin {equation}\label{26}
F(\{\l_{jk}\})=F_{j}, \quad F(B)=\sum_{\l_{jk}\in B} F_{j}.
\end{equation}
define the scalar measure $F(\cd)$ on bounded Borel sets $B\subset \bR$. For
this measure we have
\begin {equation*}
F([k,k+1))=\sum\limits_{j=1}^\infty  F_j=s, \quad k\subset\bN
\end{equation*}
and, consequently, $F(\bR)=\infty$. Moreover,
\begin {equation*}
\int\limits_{(-\infty,1)}\frac {d F(t)} {t^2+1}=0 \;\; \text{and}\;\;
\int\limits_{[k,k+1)}\frac {d F(t)} {t^2+1}=\sum_{j=1}^\infty \frac
{F_j}{\left(k+\frac 1 {j+1}\right)^2+1} \leq \sum_{j=1}^\infty
\frac{F_j}{k^2}=\frac 1 {k^2} s
\end{equation*}
for all $k\in \bN$, which implies that
\begin {equation*}
\int\limits_\bR\frac {d F(t)} {t^2+1}=\sum_{k=1}^\infty
\int\limits_{[k,k+1)}\frac {d F(t)} {t^2+1}\leq s \sum_{k=1}^\infty \frac 1
{k^2} <\infty.
\end{equation*}
Hence the measure $F(\cd)$ satisfies the relation \eqref{22}. Next for a given
$\l_k=k$ consider the interval $(\a_k,\b_k)$ such that $(k-1)+\tfrac 1 2
<\a_k<k$ and $k+ \frac 1 2<\b_k<k+1$. Then
\begin {equation*}
\int\limits_{[\a_k,\b_k)}\frac {d F(t)} {(t-\l_k)^2}=\sum_ {j=1}^\infty\frac {
F_j} {\left [\left(k+\frac 1 {j+1} \right)-k\right]^2}=\sum _{j=1}^\infty
(j+1)^2 F_j<\infty
\end{equation*}
and by \eqref{22} the inequality \eqref{25} is also valid. This gives the
relation \eqref{23} for the measure \eqref{26}.

Thus in both cases (i) and (ii) we constructed the countable infinite subsets
$Y_1=\{\l_{jk}\}$ and $Y_2=\{\l_k\}$ of $I$ and the discrete measure $F(\cd) $
with the following properties: 1) $Y_1$ consists of isolated points and $Y_2$
is the set of all limit points of $Y_1$ belonging to $I$;  2) the measure
$F(\cd)$ is concentrated on $Y_1$, $F(\bR)=\infty$ and the relations \eqref{22}
and  \eqref{23} are satisfied. This properties imply that the multiplication
operator $\wt A_F$ in $L_2(F)$ (see \eqref{12.1}) satisfies the equality
\begin {equation}\label{27}
\s_e (\wt A_F)\cap I=Y_2(=\{\l_k\}).
\end{equation}
Next assume that $A_F\subset \wt A_F$ is a simple symmetric operator in
$L_2(F)$ given by \eqref{12.2}.  Then $n_\pm (A_F)=1$ and in view of \eqref{27}
and Lemma \ref{lem6a} the set $\s_e (\wt A)\cap I(=\{\l_k\})$ is infinite for
any self-adjoint extension $\wt A \supset A_F$. Moreover, by \eqref{23} and
Lemma \ref{lem6}, 2) $Y_2 \subset \wt\rho_I (A_F)$ and the equality \eqref{13a}
gives
\begin {equation*}
I\setminus Y_2=\hat \rho_I(A_F)\subset\wt\rho_I (A_F).
\end{equation*}
This implies that $I=\wt\rho_I (A_F)$ and hence $A_F$ is a desired operator.

In the case of an arbitrary $d\in\bN$ we put $A=\bigoplus\limits_{k=1}^d A_F$,
where $A_F$ is the constructed above simple symmetric operator with $n_\pm
(A_F)=1$. It is clear that the operator $A$ has the required properties.
\end{proof}
\section{Differential operators}
In this section we apply the obtained results to  differential operators
generated by the formally self-adjoint differential expression
\begin {equation}\label{29}
l[y]=\sum_{k=1}^n (-1)^k ( (p_{n-k}y^{(k)})^{(k)}-\tfrac {i}{2}
[(q_{n-k}^*y^{(k)})^{(k-1)}+(q_{n-k} y^{(k-1)})^{(k)}])+p_n y
\end{equation}
of an even order $2n$. The coefficients $p_k(\cd)$ and $q_k(\cd)$ of this
expression are defined on an interval $(a,b), \; -\infty\leq a <b\leq\infty,$
take on values in $[\bC^m]$ and possess the following properties:

(a) $p_k,\;q_k$ are measurable on $(a,b)$;

(b) $p_k(t)=p_k^*(t)\; (k=0\div n)$ and $0\in\rho (p_0(t))$ almost everywhere
on $(a,b)$;

(c) the operator functions $p_k\; (k=2\div n),\; q_k\;(k=1\div (n-1)), \;
p_0^{-1}, \; q_0^*p_0^{-1}$ and $\tfrac 1 4 q_0^* p_0^{-1}q_0-p_1$ are locally
integrable on $(a,b)$.

The expression \eqref{29} is called regular at $a$, if $a>-\infty$ and the
assumptions on the coefficients are satisfied in $[a,b)$ instead of $(a,b)$.
The regularity of \eqref{29} at $b$ is defined correspondingly.

Next assume that $y^{[k]}(\cd), \; k=0\div 2n$ are the quasi-derivatives of a
function $y(\cd):(a,b)\to \bC^m$ \cite{Nai,KogRof74} and let $\cD (l)$ be the
set of all functions $y(\cd)$ such that  the  quasi-derivatives
$y^{[k]}(\cd),\; k=0\div (2n-1)$  are absolutely continuous  in $(a,b)$. Then
for each function $y\in\cD (l)$ the equality $l[y]=y^{[2n]}$ is valid.

For a given interval $(\a,\b)\subset \bR$ denote by $L_2(\a,\b)$ the Hilbert
space of all measurable  functions $f(\cd):(\a,\b)\to \bC^m$ such that
$\int\limits _\a^\b ||f(t)||^2\,dt<\infty$. As  is known \cite{Nai,Wei}  the
expression \eqref{29} generates the maximal operator $L$ in $L_2(a,b)$ defined
on the domain $\cD (L):=\{y\in\cD (l)\cap L_2(a,b): l[y]\in L_2(a,b)\}$ by
$Ly=l[y], \; y\in\cD (L)$. Moreover, the minimal operator $L_0$  is defined by
$L_0=\ov{L_0'}$, where $L_0'$ is a restriction of $L$ onto the linear manifold
of all functions $y\in \cD (l)$ with compact support. It is known \cite
{Nai,Wei} that $L_0$ is a  densely defined symmetric operator  in $L_2(a,b)$
and $L_0^*=L$.

 For a given point $c\in (a,b)$ denote by $l_a[y]$ and $l_b[y]$ the
restrictions of the expression $l[y]$ onto the intervals $(a,c)$ and $(c,b)$
respectively and let $L_{a0} \; (L_{b0})$ be the minimal operator in $L_2(a,c)$
(resp. $L_2(c,b)$) generated by $l_a[y]$ (resp. $l_b[y]$). It is clear that for
each $\l \in\bC$ the defect subspace $\gN_\l(L_{a0})$  $(\gN_\l(L_{b0}))$ is
the set of all solutions of the equation
\begin {equation}\label{30}
l[y]-\l y=0,
\end{equation}
which lie in $L_2(a,c)$ (resp. $L_2(c,b)$). Therefore the defect number
$n_+(L_{a0})$  $(n_+(L_{b0}))$ can be defined as the number of linearly
independent solutions of the equation \eqref{30} with $\l=i$ belonging to
$L_2(a,c)$ (resp. $L_2(c,b)$). Similarly  one defines (with $\l=-i$ in place of
$\l=i$) the defect numbers $n_-(L_{a0})$ and $n_-(L_{b0})$.

If the operators $L_{a0}$ and $L_{b0}$ have equal deficiency indices
\begin {equation}\label{30a}
n_+(L_{a0})=n_-(L_{a0})=:d_a, \quad n_+(L_{b0})=n_-(L_{b0})=:d_b,
\end{equation}
then $n\,m\leq d_a\leq 2n\,m, \; \;n\,m\leq d_b\leq 2n\,m$ and the operator
$L_0$ also has equal deficiency indices
\begin {equation*}
n_+(L_0)=n_-(L_0)=d_a+d_b-2n\,m.
\end{equation*}
In this connection note that the relations \eqref{30a} hold if $m=1$ (the
scalar case) and in formula \eqref{29} $q_k=0$. Observe also that all the above
definitions and assertions  do not depend on the choice of the point $c\in
(a,b)$.

Application of Theorem \ref{th7} to the minimal differential operator $L_0$
gives the following result.
\begin{theorem}\label{th9}
Let the operators $L_{a0}$ and $L_{b0}$ have equal deficiency indices
\eqref{30a} and let  $I=(\mu_1,\mu_2), \; -\infty\leq\mu_1 <\mu_2\leq \infty,$
be an interval such that for some (equivalently, for all) $c\in (a,b)$ and for
all $\l\in I$, besides an at most countable set $X\subset I$, the equation
\eqref{30} has $d_a$ linearly independent solutions belonging to $L_2(a,c)$ and
$d_b$ linearly independent solutions which lie in $L_2(c,b)$. Then for any
self-adjoint extension $\wt A\supset L_0$ the statement 1) of Theorem \ref{th7}
holds.
\end{theorem}
\begin{proof}
Since the expressions $l_a[y]$ and $l_b[y]$ are regular at $c$, it follows that
the corresponding minimal operators $L_{a0}$ and $L_{b0}$ are simple (see for
instance \cite{Gil72}). Hence the symmetric operator $\hat L_0:= L_{a0}\oplus
L_{b0}$ in $L_2(a,b)$ is also simple and in view of the equality $\gN_\l(\hat
L_0)=\gN_\l( L_{a0})\oplus \gN_\l( L_{b0})$ one has
\begin {equation*}
\dim\gN_\l(\hat L_0)=\dim\gN_\l( L_{a0})+\dim \gN_\l( L_{b0}), \quad \l\in\bC.
\end{equation*}
Therefore by \eqref{30a} $n_\pm (\hat L_0)=d_a+d_b$ and
\begin {equation*}
\dim\gN_\l(\hat L_0)=d_a+d_b(=n_\pm (\hat L_0)), \quad \l\in I\setminus X,
\end{equation*}
which implies that $\wt \rho_I(\hat L_0)=I\setminus X$. Moreover, $\hat
L_0\subset L_0$ and consequently $\hat L_0\subset \wt A$ for any self-adjoint
extension $\wt A\supset L_0$. Now it remains to apply Theorem \ref{th7} to
$\hat L_0$ and any self-adjoint extension $\wt A\supset L_0(\supset \hat L_0)$.
\end{proof}
The following corollary is immediate from Theorems \ref{th9} and \ref{th7}.
\begin{corollary}\label{cor10}
Let the expression \eqref{29} be regular at $a$ and let the corresponding
minimal operator $L_0$ has equal deficiency indices $d=n_\pm (L_0)$. Moreover,
let $I=(\mu_1,\mu_2), \; -\infty\leq\mu_1 <\mu_2\leq \infty,$ be an interval
such that Eq. \eqref{30} has $d$ linearly independent solutions which lie in
$L_2(a,b)$ for all $\l\in I$ besides an at most countable set $X\subset I$.
Then:

1)for any self-adjoint extension $\wt A\supset L_0$ the statement 1) of Theorem
\ref{th7} holds;

2) the set of all points $\l\in I$ such that $\ov{\cR (L_0-\l)}\neq \cR
(L_0-\l)$ is nowhere dense in $I$.
\end{corollary}
The particular case of \eqref{29} is the scalar  Sturm - Liouville expression
\begin {equation} \label{31}
l[y]=-y''+p(t)y, \quad t\in (0,\infty),
\end{equation}
where $p(\cd):(0,\infty)\to \bC$ is a scalar function such that
$p(t)=\ov{p(t)}$ and $p(t)\in L_1(0,c)$ for every $c\in (0,\infty)$ (this means
that the expression \eqref{31} is regular at 0). Let the minimal operator $L_0$
of the expression \eqref{31} has minimal deficiency indices $n_\pm (L_0)=1$
(the limit point case). For a given $\t\in\bR$ consider the boundary value
problem defined by the equation \eqref{30} and the boundary condition
\begin {equation} \label{32}
y'(0)-\t y(0)=0.
\end{equation}
Assume that $\f (t,\l)$ is the solution of \eqref{30} with the initial data $\f
(0,\l)=1,\; \f'(0,\l)=\t$ and let $\wt A_\t$ be a self-adjoint extension of
$L_0$ with the domain
\begin {equation*}
\cD (\wt A_\t)=\{y\in\cD (L): y'(0)=\t y(0)\}.
\end{equation*}
As is known \cite{Nai} the scalar measure $F(\cd):\cB_b\to \bR$ is called a
spectral measure of the boundary problem \eqref{30}, \eqref{32} if the relation
(the Fourier transform)
\begin {equation} \label{33}
L_2(0,\infty)\ni f \to (V f)(\l)=\int_0^\infty \f(t,\l)f(t)\, dt\,\in L_2(F)
\end{equation}
defines the unitary operator $V:L_2(0,\infty)\to L_2(F)(=L_2(F,\bC))$ such that
the operators $\wt A_\t$ and $\wt A_F$ (see \eqref{12.1}) are unitary
equivalent by means of $V$.

If $F(\cd)$ is a spectral measure of the problem \eqref{30}, \eqref{32}, then
the unitary operator  $V$ \eqref{33} gives the unitary equivalence between the
minimal operator $L_0$ and the symmetric operator $A_F$ defined by
\eqref{12.2}. Observe also that $F(\cd)$ is  the spectral measure of the
Titchmarsh--Weyl function $m(\cd)$ of the boundary problem \eqref{30},
\eqref{32}\cite{Nai} and hence it satisfies the relation \eqref{2}.

In the following proposition we show that the conjecture by Hartman and Wintner
on the spectrum of a self-adjoint Sturm - Liouville operator is false (for more
details see Introduction).
\begin {proposition}\label{pr11}
For any finite interval $I=(\mu_1,\mu_2), \; -\infty <\mu_1 <\mu_2 < \infty,$
there exists a Sturm--Liouville expression \eqref{31} such that the deficiency
indices of the minimal operator $L_0$ are $d=n_\pm (L_0)=1$ and the following
statements hold:

1) for all $\l\in I$ Eq. \eqref{30} has the unique solution which lies in
$L_2(0,\infty)$;

2) for any self-adjoint extension $\wt A\supset L_0$ the interval $I$ contains
infinitely many points of the essential spectrum $\s_e (\wt A)$.
\end{proposition}
\begin{proof}
Let $I=(\mu_1,\mu_2)$ be a finite interval, let $\{\l_k\}. \; \{\l_{jk}\}$ and
$\{F_{jk}\}$ be the same as in the proof of Proposition \ref{pr8} (case (i))
and let $F'(B)$ be the measure on Borel sets $B\subset I$ defined by
\eqref{20}. According to \cite[ch. 8. 26.3]{Nai} the measure $F'(\cd)$ can be
extended to the measure $F(\cd)$ on bounded Borel sets $B\subset \bR$ with the
following property: there exists the Sturm--Liouville expression  \eqref{31}
and a real $\t$ such that the corresponding minimal operator $L_0$ has the
deficiency indices  $d=n_\pm (L_0)=1$ and $F(\cd)$ is the spectral measure of
the boundary problem \eqref{30}, \eqref{32}.

Next assume that $A_F$ and $\wt A_F$ are the operators \eqref{12.1} and
\eqref{12.2} respectively. Then repeating the same reasonings  as in the proof
of Proposition \ref{pr8} one obtains that $\wt\rho_I (A_F)=I$ and for any
self-adjoint extension $\wt A\supset A_F$ the set $\s_e (\wt A)\cap I$ is
infinite. Since the minimal operator $L_0$ is unitary equivalent to $A_F$ (by
means of the Fourier transform \eqref{33}), the operator $L_0$ has the same
properties as $A_F$. This implies the desired statements 1) and 2).
\end{proof}

       \end{document}